\documentclass[12pt, reqno]{amsart}

\usepackage{fullpage}
\usepackage{amsmath}
\usepackage{amssymb}
\usepackage{dsfont}
\usepackage{tikz}
\usepackage[english]{babel}
\usetikzlibrary{calc}
\usepackage{pgfplots}
\usepackage{microtype}
\usepackage{booktabs}
\usepackage{hyperref}
\usetikzlibrary{backgrounds,datavisualization.formats.functions}

\newcommand{\Chi}{\boldsymbol{\chi}}
\newcommand{\NN}{\mathbb{N}}
\newcommand{\QQ}{\mathbb{Q}}
\newcommand{\RR}{\mathbb{R}}
\newcommand{\CC}{\mathbb{C}}
\newcommand{\EE}{\mathbb{E}}
\newcommand{\VV}{\mathbb{V}}
\newcommand{\PP}{\mathbb{P}}
\newcommand{\ZZ}{\mathbb{Z}}

\newcommand{\TT}{\mathbb{T}}

\newcommand{\mc}{\mathcal}

\newcommand{\Vol}{\text{Vol}}
\DeclareMathOperator*{\esssup}{ess\,sup}

\newtheorem{lemma}{Lemma}[section]
\newtheorem{theorem}[lemma]{Theorem}
\newtheorem{corollary}[lemma]{Corollary}

\theoremstyle{definition}

\definecolor{brickred}{rgb}{0.8, 0.25, 0.33}

\numberwithin{equation}{section}

\title{On the randomized lattice-point count problem}
\author{Steffen Schlie}
\address{Universität Osnabrück, Fakultät für Mathematik,
  Albrechtstraße 28a, 49076 Osnabrück, Germany}
\email{steffen.schlie@uni-osnabrueck.de}
\thanks{This work was supported by the German Research Foundation (DFG) via SPP 2458 \emph{Combinatorial Synergies}, project number 539849618.}

\begin{document}

\begin{abstract}
  Counting (integer) lattice points in a convex body $K$ has a rich history, beginning with the famous circle problem imposed by Gauss in 1834. Clearly, the number of covered lattice points depends on the position, orientation, size, and shape of $K$ in a delicate way. In 1947, Kendall was the first to consider a randomized version of this problem by allowing random translations of $K$. In this survey, we aim to provide an overview of the current state of research on the \emph{randomized lattice-point count problem}, i.e.\ on the distribution of the lattice-point count under random rotations and translations of $K$. As it turns out, even the variance is far from being well-understood, leaving multiple open questions for future research. The motivation for the present work grew out of the author's research interest in a related combinatorial problem, dealing with the convex hull of all covered lattice points.
\end{abstract}

\maketitle


\section{Introduction}

Suppose that you are given a $d$-dimensional compact convex set $K$ with non-empty interior. An innocent-looking but surprisingly profound questions asks for the number of (integer) lattice points lying inside $K$. In other words, what implications do the size, shape, position and orientation of $K$ impose on its discrete volume? Let us consider $d=1$ as an introductory example.

The only $1$-dimensional compact convex sets are closed intervals $[a,b]$. Little effort is required to check that there are $\lfloor b \rfloor - \lceil a \rceil +1$ lattice points inside such an interval. If we identify the position of $[a,b]$ through $c=(a+b)/2$ and slide it over the real line, the resulting lattice-point count is a step-function in $c$ with \grqq jumps\grqq\ from $\lfloor b-a\rfloor$ to $\lceil b-a\rceil$ (and vice versa) at values $c = k \pm (b-a)/2$, where $k \in \ZZ$. Due to regularity properties of $\ZZ$ we already obtain a complete answer to the above question by considering $c \in [0,1)$.

Deriving a similar characterization in the plane is much more complex and versatile. Compact convex sets of dimension $2$ can take all kinds of shapes and forms - polygons, circle or everything general - and be arbitrarily positioned and oriented in the plane. While an enormous amount of ideas and results have been contributed to first instances of the problem which date back almost two centuries, many initial conjectures remain unconquered. The most prominent one is known as \textbf{Gauss' circle problem}: What is the number of lattice points inside a planar circle of radius $\lambda \geq 0$ centered at the origin? It is easily proved that, as $\lambda$ tends to infinity, the number of interest denoted by $N(\lambda)$, is asymptotically of the area of the circle $\pi\lambda^2$. However, there is a fluctuation $\mc{D}(\lambda)$ around this value. The core of the Gauss circle problem has since been to find a correct upper bound for the modulus of this error term. Gauss \cite{Gauss1834} found a first estimate in $|\mc{D}(\lambda)| < 2\sqrt{2\pi\lambda}+2\pi$, implying $|\mc{D}(\lambda)|=\mc{O}(\lambda)$. Later, Hardy \cite{Hardy1915} proved that $|\mc{D}(\lambda)|=\mc{O}(\lambda^\alpha)$ for some value $1/2 \leq \alpha <1$, with $1/2$ being the exact lower bound. This was proved independently by Landau \cite{Landau1927}. Moreover, the conjecture associated with the Gauss circle problem goes back to Hardy again:\\

\noindent \emph{Conjecture 1:} Let $K \subset \RR^2$ be a circle of radius $\lambda > 0$ centered at the origin. Then $|\mc{D}(\lambda)| = \mc{O}(\lambda^{1/2+\varepsilon})$ for all $\varepsilon > 0$.\\

Vorono\"{\i} \cite{Voronoi1903}, Sierpi\'nski \cite{Sierpinski1906} and van der Corput \cite{vanderCorput1922} improved the upper bound to $\alpha < 2/3$. Subsequently, various intermediate refinements have been made - in chronological, but by no means complete order - by van der Corput \cite{vanderCorput1928}, Littlewood, Walfisz and Landau \cite{LittlewoodWalfisz1924}, Titchmarsh \cite{Titchmarsh1935}, Hua \cite{Hua1942}, Chen \cite{Chen1963}, Kolesnik \cite{Kolesnik1969,Kolesnik1982}, Iwaniec and Mozzochi \cite{IwaniecMozzochi1988} and Huxley \cite{Huxley1993}. The best result we have today is known to be $\alpha < 131/208$ by Huxley \cite{Huxley2003}, which is still far from proving Conjecture 1.

It didn't take long for the Gauss circle problem to be generalized to $d\geq 2$ and a wider range of compact convex sets. The study of this setting began with groundbreaking work by Minkowski \cite{Minkowski1910}, which became the basis for the branch of number theory called geometry of numbers.

\begin{theorem}[Minkowski's Lattice Point Theorem \cite{Minkowski1910}]
    Let $L$ be a lattice of determinant $\det(L)$ in $\RR^d$ and let $K \subset \RR^d$ be a convex set that is symmetric with respect to the origin. If $\operatorname{Vol}_d(K) > 2^d \det(L)$, then $K$ contains at least three lattice points.
\end{theorem}

Note that this is just an existence result for a specially shaped convex body. A lot of subsequent research in the 20th and 21st century dealt with the following: Let $K \subset \RR^d$ be a convex body (i.e.\ a compact convex set with non-empty interior) and $\lambda > 0$. What is the correct upper bound for the absolute discrepancy $|\mc{D}(\lambda)|$ between the number of integer points $N(\lambda)$ inside $\lambda K$ and its volume $\Vol_d(\lambda K)$ as $\lambda$ tends to infinity? We refer to this question as the \textbf{lattice-point count problem}.

Even simple two-dimensional examples (e.g., an axis-parallel rectangle vs.\ an ellipse) show that it is worthwhile to consider different classes of convex bodies separately. At one extreme are lattice polytopes $P \subset \RR^d$. Ehrhart's \cite{Ehrhart1962} famous theorem states that $N(n)$ is a polynomial of degree $d$ in $n \in \NN_{> 0}$. Even more, its formal power series called \emph{Ehrhart series} turns out to be a rational function $\operatorname{Ehr}_P(z) = 1 + \sum_{n \geq 1} N(n)z^n = h^\ast_P(z)/(1-z)^{d+1}$. The polynomial $h^\ast_P(z)$ has degree at most $d$ and is referred to as the $h^\ast$-polynomial of $P$. Note that there are similar results for rational polytopes. Understanding the coefficients of $h^\ast_P(z)$ is a vibrant field today and the interested reader may be drawn to \cite{Beck2010} for an introductory text on Ehrhart theory.

In general, it is not known what the geometric interpretation of the coefficients of $h^\ast_P(z)$ is. But there exists a neat formula for $d=2$ discovered by Pick \cite{Pick1899}:

\begin{theorem}[Pick's Formula \cite{Pick1899}]
    Let $P \subset \RR^2$ be a lattice polygon of area $A$ and let $n \in \NN_{> 0}$. Further, let $B$ denote the number of boundary lattice points of $P$. Then $N(n) = An^2+(B/2)n+1$.
\end{theorem}

The other extreme is represented by the class of smooth convex bodies $K$ with everywhere positive Gaussian curvature (see Section 2 for a definition). A conjecture in general dimension still open to this day is the following:\\

\noindent\emph{Conjecture 2:} Let $K \in \mc{K}^d$ be a smooth convex body. Then $|\mc{D}(\lambda)|=\mc{O}(\lambda^{d-2+\alpha(d)})$ with $\alpha(d) = 0$ for $d\geq 5$ and $\alpha(d) = \varepsilon > 0$ for $d=3$ and $4$ with an arbitrary $\varepsilon$.\\

Initial results on the aforementioned lattice-point count problem by Siegel \cite{Siegel1935}, Hlawka \cite{Hlawka1950} or Krätzel and Nowak \cite{KraetzelNowak1991,KraetzelNowak1992} led to the formulation of Conjecture 2 by M\"uller \cite{Müller1999}. Most recently, Guo \cite{Guo2012} proved the sharpest estimates to date. 

Considering a randomized setting of the lattice-point count problem is worthwhile for many reasons. For instance, aiming for generic results averages out extreme positions of orientations, even when their share among all possible configurations is quite small. The central problem in this paper is the \textbf{randomized lattice-points count problem}: As $\lambda$ tends to infinity, what is the limit distribution of the number of lattice points and its discrepancy considered as random variables with respect to a random translation and/or orientation? We will make this more precise in Section 3.

Another perspective on counting lattice points in randomly positioned convex bodies was introduced by Bárány and Matou\v{s}ek \cite{BaranyMatousek2005}, who went one step further and investigated the generic geometric and combinatorial structure described by the convex hull of all lattice points inside $K$. This perspective is of particular interest to the author and that is why we will have more to say about it in the last section.

This paper is structured as follows. In Section 2 we recall basic definitions pertaining to convex bodies and concepts from harmonic analysis. In Section 3 we review classical results on the expected value of the lattice-point count and add a series of upper bound estimates for the discrepancy. Section 4 is devoted to the variance, which might be the most studied parameter in this setting. We also adjust the setting by considering thin annuli as observation windows, which is mostly motivated by physics. Section 5 deals with higher moments and $p$-th power estimates. Section 6 discusses some of the tools and techniques used in the proofs. As the variance is not understood with exact bounds, we emphasize possible ends where further investigations could be fruitful. Finally, Section 7 summarizes applications from other sciences with a focus on a combinatorial problem which is of interest to us.

\section{Definitions}

A convex body $K$ of dimension $d$ is a compact convex subset of $\RR^d$ with non-empty interior $K^\circ$. Its boundary is denoted by $\partial K$ and we define $\mathcal{K}^d$ to be the set of all convex bodies of dimension $d$ in $\RR^d$. The $d$-dimensional volume $\Vol_d(K)$ is the classical Lebesgue measure, with $\Vol_{d-1}(K)$ measuring the surface area. At several points in Section 6 we will need the characteristic function $\Chi_K$ of a convex body $K$. We say that $x \in \partial K$ is smooth if there is a unique supporting hyperplane of $K$ at $x$. A convex body $K$ - or more precisely its boundary $\partial K$ - is called smooth if all boundary points are smooth. Equivalently, $\partial K$ is at least twice continuously differentiable. The Gaussian curvature $\kappa(x)$ is the product of all principal curvatures, i.e.\ the product of the eigenvalues of the shape operator at that point. Points with $\kappa(x) = 0$ are called flat. We say that a convex body has $C_+$-boundary whenever its boundary is smooth with everywhere positive Gaussian curvature. When we require $\partial K$ to be $k$-times differentiable, we write $C_+^k$. An object extensively studied and associated to $K$ is the annulus of thickness $r > 0$, $K(r) = (1+r/2)K \setminus(1-r/2)K$. 

A polytope $P$ is the convex hull of finitely many points $z_1,\dots,z_m$ in $\RR^d$ and hence a convex body by definition. The set of all polytopes in $\RR^d$ is denoted by $\mathcal{P}^d$. Polytopes in the plane are called polygons. We say that $P$ is a lattice polytope (rational polytope) if $z_1,\dots,z_m \in \ZZ^d \; (\QQ^d)$. The dimension of $P$ is the dimension of the affine space spanned by $P$. A face of $P$ is a set of the form $P \cap H$, where $H$ is a supporting hyperplane of $P$. The combinatorial structure of a polytope is described by the $f$-vector of $P$, $f(P) := (f_0(P),f_1(P),\dots,f_{d-1}(P))$, where $f_i(P)$ denotes the number of $i$-dimensional faces of $P$. 

The slopes and affine lengths of the sides of a polygon $P$ play a crucial role in many studies. Thus, we shall use the following definitions. For a side of a polygon $P$, we call the affine length $l$ the number of segments into which lattice points divide that side. Let $n_i = (\cos(\alpha_i),\sin(\alpha_i))$ denote the outward unit normal of the $i$-th side of $P$. Then $\alpha_i$ is the slope of that side. $P$ is called irrational if all of the slopes $\alpha_1,\dots,\alpha_h$ are irrational numbers for $j=1,\dots,h$ ($h$ is the number of sides of $P$). We call $P$ algebraic, if all slopes are irrational algebraic numbers.\\

Since many of the tools and techniques in Chapter 6 originate from Fourier analysis, we will also recall the most important definitions here. We write $L^p(\RR^d)$ for the space of measurable functions $f : \RR^d \to \CC$ for which
\begin{align*}
    \Vert f\Vert_{L^p} = \left( \int_{\RR^d} |f(x)|^p \right)^{1/p} < \infty
\end{align*}
when $1 \leq p < \infty$, and
\begin{align*}
    \Vert f \Vert_{L^\infty} = \esssup_{x \in \RR^d} |f(x)|
\end{align*}
when $p = \infty$. For $f \in L^1(\RR^d)$ its Fourier transform is defined by $\hat{f}(\xi) = \int_{\RR^d} f(x) e^{-2\pi i x\cdot \xi}dx$, with inverse transform given by $f(x) = \int_{\RR^d} \hat{f}(\xi) e^{2\pi i x \cdot \xi}d\xi$ whenever both sides are integrable. For a 1-periodic function $g \in L^1(\RR)$, its Fourier coefficients are $\hat{g}(n) = \int_0^1 g(x) e^{-2\pi i n\cdot x}dx$, and the corresponding Fourier series is $\sum_{n \in \ZZ} \hat{g}(n) e^{2\pi i n\cdot x}$.

Throughout this paper, $K$ denotes a convex body and $\lambda > 0$ a dilation factor.

\section{Basic results and discrepancy bounds}

Let $SO(d)$ denote the special orthogonal group, whose elements $\rho$ are interpreted as rotations in dimension $d$. Similarly, let $\TT^d \simeq \RR^d/\ZZ^d$ denote the unit torus, whose elements $t$ we interpret as translation vectors. The set $SO(d) \times \TT^d$ becomes a probability space if we equip it with the product measure of the normalized Haar measure for groups and the Lebesgue measure restricted to the torus. In order to randomize the location and orientation of $K$ in $\RR^d$, we uniformly choose elements from this probability space. We use uppercase letters $(R,T)$ when considering these pairs as random variables and lowercase letters $(\rho,t)$ for fixed, non-random choices from $SO(d) \times \TT^d$. 

The main objects of our study, which we define for a randomized choice $(R,T)$, are the \textbf{lattice-point enumerator}
\begin{align}\label{eq:LatticePointEnumerator}
    N_{R,T}(\lambda) := \#(R(\lambda K + T) \cap \ZZ^d)
\end{align}
and the \textbf{discrepancy}
\begin{align}\label{eq:Discrepancy}
    \mc{D}_{R,T}(\lambda) := N_{R,T}(\lambda) - \Vol_d(\lambda K).
\end{align}

Note that these are random variables with respect to the probability space defined above. We can now formulate our central question: As $\lambda \to \infty$, what is the limit distribution of $N_{R,T}(\lambda)$? 

The expected value $\EE N_{R,T}(\lambda)$ is well understood, so that the following fundamental result can be found in various references on this topic.

\begin{theorem}\label{th:ExpectedValue}
    \cite{Kendall1947} Let $K \in \mathcal{K}^d$ be arbitrary and $\lambda > 0$. Then $\EE N_{R,T}(\lambda) = \textnormal{Vol}_d(\lambda K)$ and $\EE \mathcal{D}_{R,T}(\lambda) = 0$.
\end{theorem}
In Section 6 we will show that the same holds for $N_T(\lambda)$ and $\mc{D}_T(\lambda)$. Even more, Theorem \ref{th:ExpectedValue} is true for arbitrary translation invariant random mosaics. Unfortunately, no such results hold for either $N_R(\lambda)$ or $\mc{D}_R(\lambda)$. This motivated investigations of $\mc{D}_\rho(\lambda)$ and $\mc{D}_{\rho,T}(\lambda)$ for fixed $\rho$ with the goal to find bounds that hold for a \grqq large number" of rotations. We briefly mention the most important results here without any details.

For the class $\mathcal{P}^2$, the trivial bound $|\mc{D}_\rho(\lambda)| \ll \lambda$ is best possible if at least one of the normals to the sides of a polygon $P$ is proportional to a vector of $\ZZ^2$. It turned out that for irrational polygons the discrepancy can be logarithmically small and that the appearance of logarithmically small errors on the set $\mathcal{P}^2$ is quite typical. 

\begin{theorem}\label{th:HLKhSkriganov}
    Let $P \in \mc{P}^d$ be a compact polyhedron with $d \geq 2$. For $d=2$ denote by $\alpha_1,\dots,\alpha_h$ the slopes of the sides of $P$. Then the following assertions hold:
    \begin{enumerate}
        \item[(i)] \cite{HardyLittlewood1922} If $P \in \mc{P}^2$ and if all $\alpha_1,\dots,\alpha_h$ are suitable irrational, then one has the bound $|\mc{D}_{\rho,T}(\lambda)| = \mathcal{O}(\log \lambda).$ \vspace{0.1 em}
        \item[(ii)] \cite{Khintchine1923} Let $P \in \mc{P}^2$. Then, for almost all $\alpha_1,\dots,\alpha_h \in \RR$ one has the bound $|\mc{D}_{\rho,T}(\lambda)| = \mathcal{O}\left((\log \lambda)^{1+\varepsilon}\right)$ with arbitrary small $\varepsilon>0$. \vspace{0.1 em}
        \item[(iii)] \cite{Skriganov1993} If $P \in \mc{P}^2$ and if each of the numbers $\alpha_1,\dots,\alpha_h$ is of algebraic irrationality, then one has the bound $|\mc{D}_{\rho,T}(\lambda)| = \mathcal{O}(\lambda^\varepsilon)$ with arbitrary small $\varepsilon>0$. \vspace{0.1 em}
        \item[(iv)] \cite{Skriganov1998} If $P \in \mc{P}^d$ is arbitrary, then $|\mc{D}_{\rho,T}(\lambda)| = \mathcal{O}(\log^{d-1+\varepsilon}(\lambda)$ for any $\varepsilon>0$ and for almost all rotations $\rho$.
    \end{enumerate}
\end{theorem}

The notation in Theorem \ref{th:HLKhSkriganov} should be read as follows: given a fixed polygon $P$ with slopes of sides $\alpha_1,\dots,\alpha_h$, we associate these slopes with a certain orientation of $P$ in the plane. This orientation is denoted by $\rho$. Recall that $\rho$ is fixed. The bounds by Hardy and Littlewood \cite{HardyLittlewood1922} and Khintchine \cite{Khintchine1923} were given on the base of continued fractions, while the third bound was due to a result by Schmidt \cite{Schmidt1980} as a consequence of Roth's theorem on approximations of algebraic irrationalities. 

For $K \in \mc{K}^d$ with certain smoothness or curvature properties, a first result was given by Colin de Verdière \cite{ColinDeVerdiere1977}. To state the result we assume that $K$ has sufficiently irrational normal vectors at boundary points where the curvature vanishes with order $k \geq 2$. 

\begin{theorem}\label{th:ColindeVerdiere}
    \cite{ColinDeVerdiere1977} Let $K \in \mc{K}^2$ be rotated by an angle $\rho \in \mathbb{S}^1$. Then there exists a function $\varphi(\rho) \in L^1(\mathbb{S}^1)$ such that $|\mc{D}_{\rho}(\lambda)| \leq \varphi(\rho)\lambda^{2/3}$. More precisely, for almost all rotation angles $\rho$ one has $|\mc{D}_{\rho}(\lambda)| \leq \mathcal{O}(\lambda^{2/3})$.
\end{theorem}

Independently, Tarnopolska-Weiss \cite{Tarnopolska-Weiss1978} obtained a similar result for $K$ with $\mathcal{C}^\infty$-boundary but finitely many points of zero curvature with finite order of contact. Generalizations of the results by Colin de Verdière and Tarnopolska-Weiss can be discussed in two different ways, both being described below.

\subsection{Weakened curvature assumptions} We are still considering the case $d=2$ and $K \in \mc{K}^d$ with isolated flat points. Does the estimate in Theorem \ref{th:ColindeVerdiere} hold for arbitrary convex bodies in the plane? Up to a logarithmic transgression, Brandolini et al \cite{Brandolini2001} answered this question in the affirmative.

\begin{theorem}
    \cite{Brandolini2001} Let $K \in \mc{K}^2$. Then $|\mc{D}_\rho(\lambda)| = \mathcal{O}(\lambda^{2/3}\log^{1/2 + \varepsilon}(\lambda))$ for almost every $\rho \in \mathbb{S}^1$ and any $\varepsilon>0$.
\end{theorem}

There is a huge difference between this upper bound and that for polygons derived by Skriganov \cite{Skriganov1993} in Theorem \ref{th:HLKhSkriganov}. Brandolini et al \cite{Brandolini2001} answered the question of which geometric properties are causing this gap. Furthermore, it is collectively conjectured that the exponent 2/3 is not sharp, and should be equal to $1/2 + \varepsilon$.

\subsection{Stronger curvature assumptions} Is there a convex body in the plane for which $|\mc{D}_\rho(\lambda)| = \mathcal{O}(\lambda^\alpha)$ holds for almost every rotation $\rho$ with $\alpha < 2/3$? Which kind of boundary structure and curvature behavior is required to obtain such improvement on the bounds obtained by Colin de Verdière \cite{ColinDeVerdiere1977} and Tarnopolska-Weiss \cite{Tarnopolska-Weiss1978}?

The first progress in this direction was achieved by M\"uller and Nowak, who in a series of papers studied compact planar domains bounded by closed smooth Jordan curves defined by analytic functions \cite{MüllerNowak1985,MüllerNowak1990,Nowak1984,Nowak1985}. Under suitable Diophantine conditions on the irrational slope at points of zero curvature, they obtained the desired bound with $\alpha$ depending only on the order of vanishing of the curvature. Let $\omega - 2$ denote the maximal order of vanishing of the curvature. A natural question to ask is whether there is a bound of order $\lambda^{2/3 - \zeta}$ for almost every rotation with $\zeta > 0$ that does not depend on $\omega$. 

To our own knowledge, there are only two results in this direction; both were shown by Guo \cite{Guo2012-2,Guo2012}. The first covers the two-dimensional case, the second is a generalization for any dimension $d \geq 2$.

\begin{theorem}\label{th:Guo}
    Assume $\omega \geq 3$.
    \begin{enumerate}
        \item[(i)] \cite{Guo2012-2} Let $K \in \mc{K}^2$ have $C_+$-boundary with finite order of contact. Further let $\zeta = 1/3831$. Then $|\mc{D}_\rho(\lambda)| = \mathcal{O}\left(\lambda^{2/3 - \zeta}\right)$ for almost every $\rho \in \mathbb{S}^1$. 
        \item[(ii)] \cite{Guo2012} Let $K \in \mc{K}^d$ ($d\geq 3$) with $0 \in K^\circ$. If the boundary is a $C_+$-hypersurface with finite order or contact then $|\mc{D}_\rho(\lambda)| = \mathcal{O}\left( \lambda^{d-2+2/(d+1) - \zeta_d} \right)$ for almost every $\rho \in SO(d)$, where $\zeta_d > 0$ is defined as
    \begin{align*}
        \zeta_d = \begin{cases}
            \frac{2(d-2)(d-1)d}{(d+1)\left( 6d^5+118d^4+109d^3-210d^2-119d+82 \right)} & \text{for } 3 \leq d \leq 4, \\ \frac{(d-3)(d-1)d}{2d^6+49d^5+123d^4-9d^3-167d^2-52d+30} & \text{for } d \geq 5.
        \end{cases}
    \end{align*}
    \end{enumerate}
\end{theorem}

One can verify that the bound in \cite{Guo2012-2} is better than that of Colin de Verdière \cite{ColinDeVerdiere1977} and Tarnopolska-Weiss \cite{Tarnopolska-Weiss1978}, and that it is better than the estimates on Conjecture 2 by M\"uller \cite{Müller1999} if $\omega > 427$. It is necessary to assume that the order of contact of the tangent lines at each boundary point is at least 3, since $\omega = 2$ would give the case of non-vanishing Gaussian curvature. For the higher-dimensional case, Iosevich, Sawyer and Seeger \cite{IosevichSawyerSeeger2002} proved the existence of $r > 2$ such that
\begin{align*}
    |\mc{D}_\rho(\lambda)| = \mathcal{O}\left( \lambda^{d-2+2/(d+1)}\log^{1/r}(2+\lambda) \right)
\end{align*}
for almost every $\rho \in SO(d)$. Results of this type with the same exponent $d-2+2/(d+1)$ can be found in Randol \cite{Randol19692} for convex domains with analytic boundary, and in Colin de Verdière \cite{ColinDeVerdiere1977} for general (not necessarily convex) domains if $d \leq 7$. The result by Guo \cite{Guo2012} in Theorem \ref{th:Guo} is an immediate consequence of a more involved result in that paper. Using the same method used in the proof of his results in \cite{Guo2012}, Guo was able to extend previous results from \cite{Guo2012-2} for planar convex bodies with no curvature assumptions on the boundary with an even better estimate.

\section{Variance estimates}

In this section, we build on the results of Theorem \ref{th:ExpectedValue}. Since the expected value is well understood, the next natural question is about the variance. The results presented in this section are considered exclusively for $N_T(\lambda)$ and $N_{R,T}(\lambda)$, because very little is known for rotations only (see last section). We will focus on order estimates and leave the question of more precise bounds to Section 6.

Kendall \cite{Kendall1947} was the first to investigate the variance of the lattice-point enumerator under random translations. He gave order estimates for the circle and general planar convex bodies with smooth boundary and positive curvature as well as higher-dimensional ellipses.

\begin{theorem}\label{th:Kendall}
    \cite{Kendall1947} \begin{enumerate}
        \item[(i)] If $K \in \mc{K}^2$ is the unit circle, then $\VV N_T(\lambda) = \mathcal{O}(\lambda)$ as $\lambda \to \infty$. 
        \item[(ii)] If $K \in \mc{K}^2$ has $C_+$-boundary, then $\VV N_T(\lambda) = \mathcal{O}(\lambda)$ as $\lambda \to \infty$.
        \item[(iii)] If $K \in \mc{K}^d$ is an ellipse whose axis are parallel to the lattice directions, then $\VV N_T(\lambda) = \mathcal{O}(\lambda^{d-1})$ as $\lambda \to \infty$. 
    \end{enumerate}
\end{theorem}

Kendall and Rankin \cite{KendallRankin1953} extended (i) to $d$-dimensional hyperspheres, obtaining a bound of order $\lambda^{d-1}$. Hlawka \cite{Hlawka1950} extended (ii) in the following way: 

\begin{theorem}\label{th:Hlawka}
    \cite{Hlawka1950} Let $K \in \mc{K}^d$ have $C_+$-boundary. Then $\VV N_T(\lambda) = \mathcal{O}(\lambda^{d-1})$ as $\lambda \to \infty$.
\end{theorem}

Brandolini et al \cite{Brandolini2015} recovered Hlawka's result, investigating upper $L^p(\TT^d)$ estimates for the discrepancy.

Theorem \ref{th:Kendall} states that all of the above convex bodies allow for a variance growth proportional to the surface area of the observation window. Let us flip our setting of the randomized lattice-point count problem for a moment by keeping $K$ fixed and applying random motions to the lattice $\ZZ^d$. A lattice $L$ is called hyperuniform, when
\begin{align*}
    \frac{\VV(\#(L \cap K_\lambda))}{\Vol_d(K_\lambda)} \to 0 \; \; \text{ as } \lambda \to \infty
\end{align*}
for a convex body $K$ whose size depends on $\lambda$. In the terminology of Torquato \cite{Torquato2018}, $L$ is said to be hyperuniform of class I when in fact $\VV (\#(L \cap K_\lambda)) = \mathcal{O}(\lambda^{d-1})$ as $\lambda \to \infty$. It is mentioned by Dereudre et al \cite{Dereudre2024} - and was proven by Beck \cite{Beck1987} - that for  convex bodies $K$ with $C_+$-boundary, randomly rotated and translated copies of $\ZZ^d$ are as hyperuniform as it gets, i.e.\ hyperuniform of class I, and that this is the slowest growth possible. One of the reasons is that smooth convex bodies with everywhere positive Gaussian curvature are very likely to miss the lattice axis, which leads to small changes in $N_t(\lambda)$ for small changes in $t$.

\subsection{Further variance estimates}

In order to extend the variance estimates for smooth convex bodies with everywhere positive Gaussian curvature to further classes of convex bodies, one might allow a finite set of isolated flat points. The following results are due to Brandolini et al \cite{Brandolini2020} and are split into the planar ($d=2$) and higher-dimensional case:

\begin{theorem}
    \cite{Brandolini2020} Let $K \in \mathcal{K}^2$ have $C_+$-boundary except for a finite number of isolated flat points of order at most $\gamma$.
    \begin{enumerate}
        \item[(i)] For $1 < \gamma \leq 2$ we have $\VV N_T(\lambda) = \mc{O}(\lambda)$.
        \item[(ii)] For $2 < \gamma$ we have $\VV N_T(\lambda) = \mc{O}(\lambda^{2-2/\gamma)})$.
    \end{enumerate}
\end{theorem}

\begin{theorem}\label{th:Brandolini20}
    \cite{Brandolini2020} Let $K \in \mathcal{K}^d$ have $C_+$-boundary except for a finite number of isolated flat points of order at most $\gamma$.
    \begin{enumerate}
        \item[(i)] For $1 < \gamma \leq 2$ we have $\VV N_T(\lambda) = \mc{O}(\lambda^{d-1})$.
        \item[(ii)] For $2 < \gamma$ we have $\VV N_T(\lambda) = \mc{O}(\lambda^{2(d-1)(1-1/\gamma)})$.
    \end{enumerate}
\end{theorem}

The estimate in (i) is the same as in Brandolini et al \cite{Brandolini2015}, while in the case (ii) the flat points give the main contribution.

We now add random rotations and discuss estimates on $\VV N_{R,T}(\lambda)$. Randol \cite{Randol19692} derived results for (not necessarily convex) sets in the plane with a boundary structure given as in Theorem \ref{th:Brandolini20}, showing that $\VV N_{R,T}(\lambda) = \mathcal{O}(\lambda)$. For convex bodies $K \in \mc{K}^2$ this result can also be found in Brandolini et al \cite{Brandolini2003}. We should note here that $\VV N_{R,T}(\lambda) = \mc{O}(\lambda)$ is false without integration in $T$, as the case of a disc and Hardy's $\Omega$-result (see \cite{Kraetzel1988}) show.

Although there are no higher-dimensional results for convex bodies with a finite number of isolated flat points, Gariboldi \cite{Gariboldi2020} considered the special case of a single flat point in arbitrary dimension.
\begin{corollary}\label{th:Gariboldi}
    \cite{Gariboldi2020} Let $K \in \mc{K}^d$ have $C_+$-boundary except for one flat point of order $\gamma \geq 2$. For all $\gamma \geq 2$, one has $\VV N_{R,T}(\lambda) = \mc{O}(\lambda^{d-1})$.
\end{corollary}

Is there anything known for the class $\mc{P}^d$ of polytopes/polygons? Motivated by the result of Hardy and Littlewood \cite{HardyLittlewood1922} stated in Theorem \ref{th:HLKhSkriganov}, Brandolini et al \cite{Brandolini1997} extended Theorem \ref{th:HLKhSkriganov} (i) in a probabilistic framework to higher dimensions.

\begin{theorem}\label{th:Brandolini97}
    \cite{Brandolini1997} Let $K \in \mathcal{P}^d$ be a polytope. Then $\VV N_{R,T}(\lambda) = \mc{O}\left(\lambda^{d-1}\right)$.    
\end{theorem}

It can be seen from the proofs given in \cite{Brandolini1997} that these estimates are sharp. The authors also considered domains more general than polyhedra; for instance, they obtain estimates for domains whose boundary has finite Minkowski measure. 

A recent result by Tokmachev \cite{Tokmachev2025} gave an exact and closed formula for the variance of a lattice polygon. However, it deals with the case $\lambda = 1$ and does not include any rotational average. 

\begin{theorem}
    \cite{Tokmachev2025} Assume $P \in \mc{P}^2$ has no parallel sides and $l_1,\dots,l_n$ be the affine lengths of the sides of $P$. Then
    \begin{align*}
        \VV N_T(1) = \frac{l_1^2 + \cdots + l_n^2}{12}.
    \end{align*}
\end{theorem}

\subsection{Variance estimates for thin annuli under random translations} 

Studying the behavior of (randomized) lattice-point enumerators has been a central field of research in number theory and probabilistic number theory for a long time. In the 90's, this problem got even greater importance because of some questions in physics. For instance, the investigation of the spectrum of certain quantum systems leads to the problem about the number of lattice points in a given domain. The local behavior of the spectrum, the number of eigenvalues of the Laplace-Beltrami operator in a randomly chosen interval of fixed constant length, leads to the probabilistic problem we devote this section to. For further details on the motivations from physics we refer the reader to a work by Major \cite{Major1992}, which gives a comprehensive first impression.

For a convex body $K \in \mc{K}^d$ we are interested in the distribution of $N_T(\lambda)$ in the annuli $K(r)$ of thickness $r > 0$. For the rest of this section we omit to emphasize that $N_T(\lambda)$ denotes the number of lattice points inside an annulus of thickness $r$ of a convex body $K$. Assume that $K$ is the unit circle in the plane. Computer simulations on the distribution of $N_T(\lambda)$ by Cheng and Lebowitz \cite{Cheng1991}, motivated by works of Sinai \cite{Sinai1991,Sinai1991-2}, Bleher \cite{Bleher1991} and Major \cite{Major1992}, led to the following conjecture: \\

\noindent\emph{Conjecture 3:} For large $\lambda$, the distribution of $N_T(\lambda)$ inside $K(r)$ with $r = (2c)/\lambda$, $c > 0$, $\lambda^2 > c$, is asymptotically Poissonian with parameter $4c\cdot\Vol_2(K) = 4c\pi$. \\

In a subsequent paper, Cheng et al \cite{Cheng1994} verified the above conjecture for planar convex bodies with $C_+^4$-boundary for the first two moments.

\begin{theorem}\label{th:Cheng}
    \cite{Cheng1994} Let $K \in \mc{K}^2$ have $C_+^4$-boundary and let $c > 0$ be some fixed positive number. For $\lambda > \sqrt{c}$, let $N_T(\lambda)$ denote the number of lattice points in $K(2c/\lambda)$. Then
    \begin{align}
        \EE N_T(\lambda) &= 4c\cdot\textnormal{Vol}_2(K),\label{eq:EVAnnuli}\\\lim_{\lambda \to \infty}\VV N_T(\lambda) &= 4c\cdot \textnormal{Vol}_2(K). \label{eq:VarianceAnnuli} 
    \end{align}
\end{theorem}

Even though the above theorem is formulated for convex bodies with a certain boundary structure, (\ref{eq:EVAnnuli}) holds for annuli of arbitrary convex bodies (see Section 6). Clearly, Theorem \ref{th:Cheng} is far from proving the conjecture of a Poissonian limit distribution; some deep number theoretical statement would be needed in order to prove that for every $k \in \NN$ it holds that
\begin{align*}
    \EE N_T(\lambda)(N_T(\lambda)-1)\cdots(N_T(\lambda)-k+1) = \left( 4c \cdot \Vol_2(K) \right)^k.
\end{align*}
A sketch with further details can be found in Section 4 in \cite{Cheng1994}.

Colzani et al \cite{Colzani2021} build upon Theorem \ref{th:Cheng} and investigated the variance in thin annuli in higher dimensions. 

\begin{theorem}\label{th:Colzani}
    \cite{Colzani2021} Assume that $K \in \mc{K}^d$ has $C_+$-boundary and contains in its interior the origin. Then for every $\alpha > (d-1)/(d+1)$ there exists $0 < \beta < 1$ and a positive constant $C$ such that for every $1 \leq \lambda < +\infty$ and every $0 < r \leq \lambda^{-\alpha}$ one has
    \begin{align}
        \Bigg| \int_{\TT^d} \bigg| N_T(\lambda) - \textnormal{Vol}_d\left( \lambda K(r) \right) \bigg|^2 dt - \textnormal{Vol}_d\left((\lambda K(r)\right)\Bigg| \leq C\cdot \textnormal{Vol}_d\left( \lambda K(r) \right) r^\beta. \label{eq:Colzani}
    \end{align}
\end{theorem}

There are a few things we want to remark. First, verify that (\ref{eq:Colzani}) indeed is an estimate for the variance. Second, it does not suffice to only assume that $r$ converges to zero. One has to require a suitable speed of convergence in order to avoid any of the critical regimes developed by Parnovski and Sidorova \cite{Parnovski2010}. More precisely, a dimension is called critical if both $L^1$ and $L^2$ norms of the discrepancy do not have upper and lower bounds of the same order as $\lambda \to \infty$. There are four essentially different regimes of the behavior of $r$ described in Section 1.2 in \cite{Parnovski2001}. And at last, it is not known whether the assumption $0 < r \leq \lambda^{-\alpha}$ with $\alpha > (d-1)/(d+1)$ can be weakened. 

However, it is equally unclear if the distribution of $N_T(\lambda)$ in thin annuli of convex bodies with $C_+$-boundary follows a certain limit distribution.

\section{Higher moments and $p$-th power discrepancy estimates}

In this section we summarize results on $p$-th power discrepancy estimates. We will only do so very briefly, because since the variance is not fully understood in large parts, higher moments suffer from the same problem. As mentioned before, we will try to get exact estimates for the variance in Section 6.

Recall that Kendall \cite{Kendall1947} proved $\mc{O}(\lambda)$ to be the order of the upper bound of $\VV N_T(\lambda)$ when $K \in \mc{K}^2$ has $C_+$-boundary. By definition of the discrepancy and the $L^p$-norm, one can restate Kendall's result as follows:
\begin{align}\label{eq:L2Discrepancy}
    \Vert \mc{D}_T(\lambda)\Vert_{L^2(\TT^2)} = \mathcal{O}\left( \lambda^{1/2} \right).
\end{align}

We want to address the follow-up question whether this bound extends to higher powers $p \geq 2$. Huxley \cite{Huxley2014} showed that for $p=4$ this holds up to a logarithmic factor:

\begin{theorem}\label{th:Huxley}
    \cite{Huxley2014} Let $K \in \mc{K}^2$ be as in Theorem \ref{th:Kendall} (ii). Then, as $\lambda \to \infty$, 
    \begin{align*}
        \Vert \mc{D}_T(\lambda)\Vert_{L^4(\TT^2)} = \mathcal{O}\left( \lambda^{1/2}\log^{1/4}(\lambda) \right).
    \end{align*}
\end{theorem}

Higher-dimensional results for smooth convex bodies with everywhere positive Gaussian curvature (except for finitely many isolated flat points) are summarized in the next theorem.

\begin{theorem}\label{th:BrandoliniTwice}
    Let $K \in \mc{K}^d$ have $C_+$-boundary.
    \begin{enumerate}
        \item[(i)] \cite{Brandolini2015} Let $\lambda > 0$.
        \begin{enumerate}
            \item[(1)] If $p < 2d/(d-1)$, then $\Vert \mc{D}_T(\lambda)\Vert_{L^p(\TT^d)} \leq C \lambda^{(d-1)/2}.$
            \item[(2)] If $p = 2d/(d-1)$, then $\Vert \mc{D}_T(\lambda)\Vert_{L^p(\TT^d)} \leq C \lambda^{(d-1)/2} \log^{(d-1)/(2d)}(\lambda).$
            \item[(3)] If $p > 2d/(d-1)$, then $\Vert \mc{D}_T(\lambda)\Vert_{L^p(\TT^d)} \leq C \lambda^{d(pd-p-d+1)/p(d+1)}.$
        \end{enumerate}
        \vspace{0.5em}
        \item[(ii)] \cite{Brandolini2020} Let $\lambda > 0$ and assume that $K$ has finitely many isolated flat points of order at most $\gamma$.
        \begin{enumerate}
            \item[(1)] For $1<\gamma\leq2$ we have same estimates as in (i) above.
            \item[(2)] For $2 < \gamma \leq d+1$ we have 
            \begin{align*}
                \Vert \mc{D}_T(\lambda)\Vert_{L^p(\TT^d)} = \begin{cases}
                    \mc{O}(\lambda^{(d-1)(1-1/\gamma)}) & \text{ if }1 \leq p \leq 2d/(d+1-\gamma), \\
                    \mc{O}(\lambda^{\frac{d(d-1)}{d+1}(1-2/(\gamma p))} & \text{ if } p > 2d/(d+1-\gamma).
                \end{cases}
            \end{align*}
            \item[(3)] For $\gamma > d+1$ and every $p \geq 1$ we have 
            \begin{align*}
                \Vert \mc{D}_T(\lambda)\Vert_{L^p(\TT^d)} = \mathcal{O}(\lambda^{(d-1)(1-1/\gamma)}).
            \end{align*}
        \end{enumerate}
    \end{enumerate}
\end{theorem}

The estimates in (i) were obtained by Brandolini et al \cite{Brandolini2015} and are essentially sharp for $p < 2d/(d-1)$, which is part of Theorem 5 in \cite{Brandolini2015}. The estimates in (ii) are due to Brandolini et al \cite{Brandolini2020}. For results including random rotations for smooth convex bodies with a single flat point we refer to Gariboldi \cite{Gariboldi2020} (in particular Theorems 2 and 5, where the latter proves a lower bound on the $L^p$-discrepancy). Within the setting of counting lattice points in thin annuli, Wu \cite{Wu2022} states upper bounds for $p\geq 2$ in the planar case, which to our knowledge is the only result in this direction.

\section{Tools and techniques}

In this section we turn to the rigorous proofs of some of the main results stated earlier. To prepare the ground, we begin by discussing the expected value; in doing so, we emphasize in which sense rotations are more difficult to handle. Building on these foundations, we continue by assembling several auxiliary tools that are standard in this area. In particular, we highlight the central role played by Fourier-analytic techniques, which will serve as our primary vehicle for estimating relevant quantities and translating structural information into computable forms. Using these tools, we then delve deeper into the analytic framework, with special attention to the behavior of the variance.

The first question we want to answer is the following: Why does Theorem \ref{th:ExpectedValue} only hold for $T$ or $(R,T)$, but not for $R$? Formally, we define the \textbf{expected value} of $N_{R,T}(\lambda)$ to be
\begin{align}\label{eq:ExpectedValue}
    \EE N_{R,T}(\lambda) := \int_{SO(d)} \int_{\TT^d} N_{R,T}(\lambda)dtd\rho.
\end{align}
In order to prove Theorem \ref{th:ExpectedValue}, we first focus on the inner integral over the torus. Recall that $\Chi_K$ denotes the characteristic function of $K$.
\begin{proof}[Proof of Theorem \ref{th:ExpectedValue}]
    As $\TT^d$ is equipped with the Lebesgue measure, one has
    \begin{equation} \label{eq:ExpectedValueFort}
    \begin{split}
      \EE N_T(\lambda) & = \EE \sum_{m \in \ZZ^d} \Chi_{\lambda K + T}(m) = \sum_{m \in \ZZ^d} \PP(m-T \in \lambda K) \\ &= \sum_{m \in \ZZ^d} \Vol_d\left((m-[0,1)^d)\cap \lambda K\right) = \Vol_d(\lambda K).
    \end{split}
    \end{equation}
    Since the volume is rotation invariant, integrating over $SO(d)$ proves the statement. The part for the discrepancy follows by definition.
\end{proof}

One readily sees where the difficulty with rotations lies. The probability that $R^{-1}(m)$ lies in $\lambda K$ is equal to $\Vol_{d-1}(|m|\mathbb{S}^{d-1} \cap \lambda K)$, as one can visualize quite easily for a planar convex body. Summation over all lattice points $m \in \ZZ^d$ does not yield the volume of $\lambda K$ but a rather complicated sum of ($d-1$)-dimensional intersections. An exact result as in Theorem \ref{th:ExpectedValue} is not known today.

The bridge between the discrete world of $N_{R,T}(\lambda)$ (or $N_T(\lambda)$) and the world of harmonic analysis is build by the following two well-known identities stated in almost every fundamental text on harmonic analysis:
\begin{theorem}[Poisson Summation Formula \cite{SteinWeiss1971}]\label{th:PoissonSummationFormula}
    Suppose $f \in L^1(\RR^d)$. Then the series $\sum_{m \in \ZZ^d} f(x+m)$ converges in the norm of $L^1(\TT^d)$ and has the Fourier expansion
    \begin{align*}
        \sum_{m \in \ZZ^d} \hat{f}(m) e^{2\pi i x \cdot m}.
    \end{align*}
\end{theorem}
In particular, for $x = 0$ we get $\sum_{m \in \ZZ^d} f(m) = \sum_{m \in \ZZ^d} \hat{f}(m)$. Roughly speaking, we obtain a direct correspondence between the evaluation of $f$ at lattice points and the Fourier coefficients. A similar connection for periodic functions is given by the following theorem:
\begin{theorem}[Parseval-Plancherel Identity \cite{Titchmarsh1939}]\label{th:ParsevalPlancherelIdentity}
    Let $f \in L^2(\TT^d)$ be 1-periodic. Then
    \begin{align*}
        \int_{\TT^d} |f(x)|^2 dx = \sum_{m \in \ZZ^d} |\hat{f}(m)|^2.
    \end{align*}
\end{theorem}

We now focus on the results in Section 4, present selected proofs and always keep in mind the question of whether the $\mc{O}$-bounds are even exact and whether we can possibly specify the constants explicitly. This would be desirable in order to address the question of the extent to which central limit theorems might exist.

It is almost mandatory that we start by proving Theorem \ref{th:Kendall} (i), because Kendall \cite{Kendall1947} first introduced the connection between Fourier series and discrete lattice point enumerators. In the following, the \textbf{variance} is always defined as
\begin{align}
    \VV N_{R,T}(\lambda) := \int_{SO(d)} \int_{\TT^d} \left( \mc{D}_{R,T}(\lambda) \right)^2 dtd\rho.
\end{align}

By Theorem \ref{th:PoissonSummationFormula}, the lattice-point enumerator $N_T(\lambda)$ has the Fourier series
\begin{align}\label{eq:FourierSeriesOfN}
    \sum_{k \in \ZZ^d} \widehat{\Chi}_{\lambda K}(k) e^{2\pi i k \cdot T},
\end{align}
with $k \cdot T$ denoting the scalar product and $\widehat{\Chi}_{\lambda K} (k) = \int_{\TT^d} N_T(\lambda)e^{-2\pi i k \cdot T}dt$. Note that $k=0$ gives 
\begin{align}\label{eq:ZerothCoefficient}
    \widehat{\Chi}_{\lambda K}(0) = \EE N_T(\lambda) = \Vol_d(\lambda K).
\end{align} 
With this at hand we can prove Theorem \ref{th:Kendall} (i).

\begin{proof}[Proof of Theorem \ref{th:Kendall} (i)]
    Recall that $K \in \mc{K}^2$ is the unit circle. Working with the expression of $\widehat{\Chi}_{\lambda K}$ above and the well-known fact from harmonic analysis that the Fourier transform of the characteristic function of a ball is a Bessel function $J_1(\cdot)$, one arrives at
    \begin{align}\label{eq:CoefficientsAsBessel}
        \widehat{\Chi}_{\lambda K}(k) = \lambda \frac{J_1(2 \pi \lambda \Vert k \Vert_2)}{\Vert k \Vert_2}. 
    \end{align}
    Let $r(l)$ denote the number of representations of $l \in \NN$ as the sum of two squares (see, e.g.\ \cite{Titchmarsh1924}). Applying Theorem \ref{th:ParsevalPlancherelIdentity} and (\ref{eq:ZerothCoefficient}) gives
    \begin{align}\label{eq:VarianceKendallProof}
        \VV N_T(\lambda) = \int_{[0,1)^2} \big| N_T(\lambda) \big|^2 dt - \pi^2 \lambda^4 = \lambda^2 \sum_{l=1}^{\infty} \frac{r(l)}{l} J_1^2(2\pi \sqrt{l} \lambda).
    \end{align}
    The series is absolutely convergent because $|J_1(z)| < cz^{-1/2}$ for $z >0$, and on using this for the Bessel function in (\ref{eq:VarianceKendallProof}) it follows that
    \begin{align*}
        \VV N_T(\lambda) = \mc{O}(\lambda) \text{ as } \lambda \to \infty.
    \end{align*}
\end{proof}

The availability of the Poisson Summation Formula for random translations is one of the key ingredients, allowing to transfer the discrete problem of counting lattice-points into the world of harmonic analysis. We think at this point it is appropriate to emphasize that an analogue for random rotations is not (yet) available. This is one reason why rotations remain out of reach for the moment.

Kendall \cite{Kendall1947} showed even more than what is stated in Theorem \ref{th:Kendall} (i): there exists a sequence $\lambda_n \to \infty$ for which $\VV N_T(\lambda_n) \geq c \lambda_n$, i.e.
\begin{align}
    \VV N_T(\lambda_n) = \Omega(\lambda_n) \text{ as } n \to \infty.
\end{align}
This holds since $J_1^2(z) < cz^{-1}$ and it not only decays like $z^{-1}$ but also attains values of this order infinitely often despite having zeros. Combined with the point-wise inequality 
\begin{align*}
    \VV N_T(\lambda) \geq 4\lambda^2 J_1^2(2\pi \lambda)
\end{align*}
following from the contribution of $k=(1,0)$ and $k=(0,1)$, we get a sequence of $\lambda$ along which $J_1$ does not vanish and therefore produces a lower bound of the desired order.

It would be of great interest whether there are results dealing with the sum of two Bessel functions of first kind and first order, stating that this sum can't become zero due to a (potential) linear independence of the square roots of the integers. \\

A natural next step is to examine the behavior of $\lambda^{-1} \VV N_T(\lambda)$ for large $\lambda$. Does this quantity alternate between two values or is it converging towards a single value? Kendall \cite{Kendall1947} used the asymptotic Bessel-function expansion of the form
\begin{align*}
    \VV N_T(\lambda) = \frac{\lambda}{\pi^2} \sum_{l \geq 1} \frac{r(l)}{l^{3/2}} \cos\left(2\pi \lambda \sqrt{l} - \frac{3\pi}{4}\right) + \mc{O}(1).
\end{align*}
Rewriting this leads to an expression of $\lambda^{-1}\VV N_T(\lambda)$ as an explicit \textit{almost periodic} oscillatory series with remainder $\mathcal{O}(\lambda^{-1})$. Define a constant $a$ as
\begin{align*}
    a^2 := \frac{1}{2\pi} \sum_{l \geq 1} \frac{r(l)}{l^{3/2}}.
\end{align*}
Then, as $\lambda$ tends to infinity, the upper and lower limits of indetermination of $\lambda^{-1}\VV N_T(\lambda)$ are equal to those of
\begin{align*}
    \frac{1}{2\pi^2} \sum_{l \geq 1} \frac{r(l)}{l^{3/2}} \left( 1- \sin(4\pi \lambda \sqrt{l}) \right)
\end{align*}
which clearly lies between 0 and $2a^2$. For a computation of the constant $a$ we refer to the appendix in \cite{Kendall1947}. \\

In the second part of this section, we want to extend the discussion on estimating $|\widehat{\Chi}_{\lambda K}(\xi)|$ for more general convex bodies with smooth boundary and everywhere positive Gaussian curvature in arbitrary dimension. In doing so, we hope to provide insight into how this differs from spherical bodies and is, so to speak, more difficult to achieve. The proof of Theorem \ref{th:Colzani} offers a nice structure which allows the reader to get a feeling of the overall steps used in many other works on this topic. We start with the example of the $d$-dimensional sphere $\{ |x| \leq \lambda\}$ of radius $\lambda$ discussed by Colzani et al \cite{Colzani2021}.

For $d\geq 2$, the Fourier transform of this sphere is a Bessel function of $d/2$-th order (see \cite{SteinWeiss1971}). More precisely, let $\xi \in \RR^d$, then
\begin{align*}
    \widehat{\Chi}_{\{|x|\leq \lambda\}}(\xi) = \lambda^d \widehat{\Chi}_{\{|x|\leq 1\}}(\lambda\xi) = \lambda^d |\lambda \xi|^{-d/2} J_{d/2}(2\pi \lambda |\xi|).
\end{align*}

Recall that we consider the setting of counting lattice points inside an annulus of thickness $r$. Its Fourier transform is
\begin{align}
    &\widehat{\Chi}_{\{\lambda -r/2 <|x|\leq \lambda+r/2 \}}(\xi) \nonumber \\ &= \widehat{\Chi}_{\{|x|\leq \lambda+r/2 \}}(\xi) - \widehat{\Chi}_{\{\lambda -r/2 <|x|\}}(\xi) \nonumber\\ &= \lambda^{d/2} |\xi|^{-d/2} \left( J_{d/2}(2\pi(\lambda+r/2)|\xi|) - J_{d/2}(2\pi(\lambda-r/2)|\xi|) \right) \label{eq:FTAnnulusSphere1}\\ &+ \left( (\lambda+r/2)^{d/2} - \lambda^{d/2} \right)|\xi|^{-d/2}J_{d/2}(2\pi(\lambda+r/2)|\xi|) \nonumber\\ &- \left( (\lambda-r/2)^{d/2} - \lambda^{d/2} \right)|\xi|^{-d/2}J_{d/2}(2\pi(\lambda-r/2)|\xi|).\nonumber
\end{align}
By the asymptotic expansions of Bessel functions for positive, real arguments (see \cite{Colzani2021}), we can express (\ref{eq:FTAnnulusSphere1}) as
\begin{align*}
    &\widehat{\Chi}_{\{\lambda -r/2 <|x|\leq \lambda+r/2 \}}(\xi) \\ &= 2\pi^{-1} \lambda^{(d-1)/2} |\xi|^{-(d+1)/2} \cos(2\pi \lambda |\xi| - \pi(d-1)/4)\sin(\pi r |\xi|) \\ &+ \mc{O}\left( \lambda^{(d-3)/2} r|\xi|^{-(d+1)/2} \right).
\end{align*}

Now assume that $K \in \mc{K}^d$ is smooth with everywhere positive Gaussian curvature. Its Fourier transform and that of its annulus behave similar like for a sphere, although a bit more complicated. We describe the idea to prove Theorem \ref{th:Colzani}, not without stating classical results who play a central role. First, we need the asymptotic expansion of the Fourier transform of $K$:

\begin{lemma}\label{lem:FTConvexBody}
    \cite{Stein1993} The Fourier transform of a characteristic function of a convex body $K \in \mc{K}^d$ with smooth boundary and everywhere positive Gaussian curvature for $|\xi| \to \infty$ has the asymptotic expansion
    \begin{align*}
        \widehat{\Chi}_K(\xi) = a(\xi) |\xi|^{-(d+1)/2} + E(\xi).
    \end{align*}
    If $\sigma(\pm\xi)$ are the points of $\partial K$ with outward unit normals $\pm\xi/|\xi|$, and $\kappa(\sigma(\pm\xi))$ are the Gaussian curvatures at those points, then
    \begin{align*}
        a(\xi) = &(2\pi i)^{-1} \exp(-2\pi i \sigma(-\xi)\cdot \xi - \pi i (d-1)/4) \kappa(\sigma(-\xi))^{-1/2} \\ &- (2\pi i)^{-1} \exp(-2\pi i \sigma(\xi)\cdot \xi + \pi i (d-1)/4) \kappa(\sigma(\xi))^{-1/2}.
    \end{align*}
    The remainder $E(\xi)$ satisfies the estimates
    \begin{align*}
        |E(\xi)| + |\nabla E(\xi)| \leq C|\xi|^{-(d+3)/2}.
    \end{align*}
\end{lemma}

This result is used to prove the decomposition of the Fourier transform of an annulus of thickness $r$.

\begin{lemma}\label{lem:FTAnnulus}
    \cite{Colzani2021} The Fourier transform of the annulus of thickness $r$ can be decomposed into
    \begin{align*}
        \widehat{\Chi}_{\lambda K(r)} (\xi) = A(\lambda,r,\xi) + B(\lambda,r,\xi).
    \end{align*}
    The main term is
    \begin{align*}
        A(\lambda,r,\xi) = &-\pi^{-1} \lambda^{(d-1)/2}|\xi|^{-(d+1)/2} \exp(-2\pi i \lambda \sigma(-\xi) \cdot \xi - \pi i (d-1)/4) \\ &\times \kappa(\sigma(-\xi))^{-1/2} \sin(\pi r \sigma(-\xi)\cdot \xi) \\ &+ \pi^{-1} \lambda^{(d-1)/2}|\xi|^{-(d+1)/2} \exp(-2\pi i \lambda \sigma(\xi) \cdot \xi + \pi i (d-1)/4) \\ &\times \kappa(\sigma(\xi))^{-1/2} \sin(\pi r \sigma(\xi)\cdot \xi).
    \end{align*}
    The remainder has the property that there exists $C > 0$ such that for every $\lambda |\xi| \geq 1$ and for every $0 <  r\leq \lambda$,
    \begin{align*}
        |B(\lambda,r,\xi)| \leq C \lambda^{(d-3)/2} r |\xi|^{-(d+1)/2}.
    \end{align*}
\end{lemma}

Since we want to present the main ideas only, the reader may be referred to \cite{Colzani2021} for further details. We conclude our sketch of the proof of Theorem \ref{th:Colzani} with the decomposition of the variance in Fourier space. Our aim is to emphasize the interplay between the \grqq Euclidean world" and the world of \grqq Fourier analysis" provided by the well-known Theorems \ref{th:PoissonSummationFormula} and \ref{th:ParsevalPlancherelIdentity}.

\begin{lemma}
    \cite{Colzani2021} The variance of the number of lattice points inside a shifted annulus $\lambda K(r)$ can be decomposed into
    \begin{align*}
        \int_{\TT^d} \Bigg| N_T(\lambda) - \Vol_d(\lambda K(r)) \Bigg|^2 dt = X(\lambda,r) + Y(\lambda,r) + Z(\lambda,r),
    \end{align*}
    where
    \begin{align*}
        X(\lambda,r) &= 2\pi^{-2} \lambda^{d-1} \sum_{n \in \ZZ^d \setminus \{0\}} \kappa(\sigma(n))^{-1} \sin^2(\pi r \sigma(n)\cdot n)|n|^{-d-1}, \\ Y(\lambda,r) &= -2\pi^{-2} \lambda^{d-1} \sum_{n \in \ZZ^d \setminus \{0\}} \cos(2\pi \lambda(\sigma(n) - \sigma(-n))\cdot n - \pi (d-1)/2 \\ &\times \kappa(\sigma(n))^{-1/2} \kappa(\sigma(-n))^{-1/2} \sin(\pi r \sigma(n)\cdot n) \sin(\pi r \sigma(-n)\cdot n)|n|^{-d-1}.
    \end{align*}
    The remainder $Z(\lambda,r)$ has the property that there exists a constant $C$ such if $\lambda \geq 1$ and $r \leq \lambda$ then
    \begin{align*}
        |Z(\lambda,r)| \leq C \Vol_d(\lambda K(r)| \lambda^{-1} r \log(2+1/r).
    \end{align*}
\end{lemma}

\begin{proof}
    By Theorems \ref{th:PoissonSummationFormula}, \ref{th:ParsevalPlancherelIdentity} and Lemma \ref{lem:FTAnnulus}, the variance equals
    \begin{align*}
        \sum_{n \in \ZZ^d\setminus \{0\}} &|\widehat{\Chi}_{\lambda K(r)}(n)|^2 \\ = &\sum_{n \in \ZZ^d\setminus \{0\}} A(\lambda,r,n)\overline{A(\lambda,r,n)} + \sum_{n \in \ZZ^d\setminus \{0\}} A(\lambda,r,n)\overline{B(\lambda,r,n)} \\ &+ \sum_{n \in \ZZ^d\setminus \{0\}} B(\lambda,r,n)\overline{A(\lambda,r,n)} + \sum_{n \in \ZZ^d\setminus \{0\}} B(\lambda,r,n)\overline{B(\lambda,r,n)}.
    \end{align*}
    Since $c|n| \leq \sigma(n) \cdot n \leq C|n|$ for some $C \geq c > 0$, Lemma \ref{lem:FTAnnulus} implies that
    \begin{align*}
        |A(\lambda,r,n)| &\leq C \lambda^{(d-1)/2} \min\{1,r|n|\}, \\ |B(\lambda,r,n)| &\leq C\lambda^{(d-3)/2} r |n|^{-(d+1)/2}.
    \end{align*}
    These estimates give
    \begin{align*}
        \sum_{n \in \ZZ^d\setminus \{0\}} &|A(\lambda,r,n)||B(\lambda,r,n)| \\ &\leq C\lambda^{d-2}r^2 \sum_{0 < |n| \leq 1/r} |n|^{-d} + C\lambda^{d-2}r \sum_{1/r < |n| < \infty} |n|^{-d-1} \\ &\leq C\lambda^{d-2}r^2 \log(2+1/r),
    \end{align*}
    and 
    \begin{align*}
        \sum_{n \in \ZZ^d\setminus \{0\}} |B(\lambda,r,n)|^2 \leq C\lambda^{d-3} r^2 \sum_{n \in \ZZ^d \setminus \{0\}} |n|^{-d-1} \leq C\lambda^{d-3}r^2.
    \end{align*}
    The main term is $\sum_{n \in \ZZ^d \setminus \{0\}} |A(\lambda,r,n)|^2$, and one can check that it is equal to $X(\lambda,r) + Y(\lambda,r)$.
\end{proof}

Even though Theorem \ref{th:Colzani} and its Lemmas above deal with thin annuli instead of convex bodies as a whole, one can get a feeling for the difficulties and complexities of variance estimates. It seems like explicit constants for lower and upper bound estimates could only be derived by very sophisticated analytic methods for specific convex bodies; see for instance the circle dealt with by Kendall \cite{Kendall1947}. We tried to highlight open ends and questions in the field of the randomized lattice-point count problem by discussing the Theorems \ref{th:Kendall} and \ref{th:Colzani} in a more detailed manner.

\section{Outlook}

In the discussion of the previous chapter, we saw that even the second moment (variance) of the lattice-point count is not yet understood. From a stochastic perspective, however, precise order estimates for the bounds of the variance are not only intrinsically interesting, but also a first step towards potential central limit theorems with exact convergence rates. For example, Kendall's study of the circle has shown that a deeper understanding of Bessel functions would be a starting point for finding concrete constants in the upper and lower bounds of order $\lambda$. This work tries to reveal such problems not only to us but also to the reader.

Sitting at the crossroads of mathematical areas like geometry of numbers, optimization, combinatorics and cryptography, the (randomized) lattice-point count problem is an interesting question in its own right. But it is equally motivating that it also has striking practical applications! In \emph{integer programming}, the problem can be expressed in the canonical form $\max\{ c^\top x : Ax\leq b, x \in \ZZ^d \}$, where the feasible region is a convex polyhedron $P = \{ x \in \RR^d : Ax\leq b\}$. Counting or estimating the number of integer points in $P$ gives valuable information: if the count is zero, infeasibility is certified; if it is small, then branch-and-bound or cutting plane methods can prune aggressively; if large, different strategies may be more appropriate. However, exact counting is generally NP-hard. Barvinok’s algorithm \cite{Barvinok1994} shows that for fixed dimension $d$, one can compute $\#(P \cap \ZZ^d)$ in polynomial time using short rational generating functions, but for higher dimensions one usually seeks approximations. Randomized counting heuristics could guide branching, cut selection, and feasibility checks in practical solvers, complementing classical structural results (total unimodularity, proximity theorems) in integer programming (see e.g.\ \cite{NemhauserWolsey1999,Schrijver1986}). 

Of course, there are many other applications that deal with numbers of integer points in convex domains - computer networking \cite{BalciogluWood2003,BaldoniSilvaLoeraVergne2004}, magic squares \cite{AhmedLoeraHemmecke2003}, matchings on graphs \cite{LovaszPlummer1986}, linear extensions of posets \cite{Stanley1986} or $t$-designs \cite{Moura1996}, to only name a few. For further applications in various areas of applied and pure mathematics see e.g.\ \cite{LoeraHawsHemmeckeHugginsSturmfelsYoshida2004,DiaconisGangolli1995,FienbergMakovMeyerSteele2001,FultonHarris1991,LiNgPurcell2022,MillerSturmfels2004,Regev2005}. It is worth noting that most of them actually deal with the classical lattice-point count problem without any random components. Problems in physic, however, include at least some \grqq randomizing aspect", as we have seen in Section 4. Other interesting problems can be found in the context of quantum field theory \cite{Elizalde1994}, quantum billiards \cite{Ivrii1980} and quantum chaos \cite{Marklof2003}.

There is another level from which we want to look at lattice points in randomly motioned convex bodies. Rather than simply counting them, we now investigate the geometric structure they determine. In particular, by taking their convex hull we arrive at a lattice polytope. This lattice-generated polytope is called \emph{randomized integer convex hull}. Formally, let $L=R(\ZZ^d+T)$ be a random copy of the $d$-dimensional integer lattice and take $K \subset \RR^d$ as given. The randomized integer convex hull of $K$ is defined as 
\[
K_L := \operatorname{conv}(K \cap L).
\]
Motivated and first investigated by Bárány and Matou\v{s}ek the randomized integer convex hull is a not well understood random lattice polytope. Bárány and Matou\v{s}ek proved asymptotic bounds for the expected number of vertices $\EE f_0(K_L)$ in any dimension as well as the expected missed area $\EE \operatorname{Area}(K \setminus K_L)$ in the planar case. All other entries of the expected $f$-vector $f(K_L)$ and higher moments like variances are completely unconquered until today. Understanding the variance of $N_{R,T}(\lambda)$ might help our understanding of the variance of $f_0(K_L)$, for instance. The only geometric result known is on the expected width, studied by Ngoc and Reitzner \cite{NgocReitzner2021}. It is of main research interest to understand $K_L$ and the entries of $f(K_L)$.

\bibliographystyle{plain}
\bibliography{references.bib}

\end{document}